# Analysis of Zeta Functions, Multiple Zeta Values, and Related Integrals

David H. Wohl

BS"D

October 25, 2001
[expansion of letter dated September 26, 2001]

Abstract

In this work, we begin to uncover the architecture of the general family of zeta functions and multiple zeta values as they appear in the theory of integrable systems and conformal field theory. One of the key steps in this process is to recognize the roles that zeta functions play in various arenas using transform methods. Other logical connections are provided by the the appearance of the Drinfeld associator, Hopf algebras, and techniques of conformal field theory and braid groups. These recurring themes are subtly linked in a vast scheme of a logically woven tapestry.

An immediate application of this framework is to provide an answer to a question of Kontsevich regarding the appearance of Drinfeld type integrals and in particular, multiple zeta values in:

    a) Drinfeld's work on the KZ equation and the associator;
    b) Etingof-Kazhdan's quantization of Poisson-Lie algebras;
    c) Tamarkin's proof of formality theorems;
    d) Kontsevich's quantization of Poisson manifolds.

Combinatorial arguments relating Feynman diagrams to Selberg integrals, multiple zeta values, and finally Poisson manifolds provide an additional step in this framework.

Along the way, we provide additional insight into the various papers and theorems mentioned above.

This paper represents an overall introduction to work currently in progress. More details to follow.

Contents:

Introduction.

Summary.

Acknowledgements.

---

Introduction.

A recent theorem of Terasoma [38](AG/9908045) relates Selberg's integral to MultiZeta Values by expressing a general form of this integral as an element of the homogeneous ring of MZV's. Terasoma's proof invokes use of the Drinfeld Associator, which had been previously been shown to be expressable in terms of MZV's by Le-Murakami [27].

In [30](9203002), Mulase points out that every Jacobian variety can be realized as a finite dimensional orbit of KP flows on an infinite dimensional Grassmanian, which in turn is equivalent to the claim that the KP equations characterize the Riemann theta functions associated with Jacobian varieties (E.Arabello-C.DeConcini [1]).

A theorem of Mulase [31] (MathPhys/9811023), based upon results of Kakei [23] (solv-int/9909023) and Goulden-Harer-Jackson [18] (AG/9902044) expresses the asymptotic expansion of an Hermitian Matrix integral in one of three ways:

    (1) as a Feynman diagram expansion;
    (2) using classical asymptotic analysis of orthogonal polynomials in
      the Penner Model resulting in Zeta Functions;
or   (3) as a Riemann theta function or a Selberg integral representing a
      tau-function solution of the KP equations.

It should also be noted that there are various derivations of hypergeometric
functions as solutions of KP equations – see for example Orlov-Scherbin [32].

Our new theorem is that by applying Terasoma's result, Mulase's
methods (2) and (3) are essentially equivalent. That is, the Selberg integral
of (3) can be related to the Zeta functions of (2). In addition, results
of Broadhurst-Kreimer (see for example, Hep-th/9609128 [5]) relate Feynman
diagrams in terms of MZV's as well. Therefore the 3 methods listed
by Mulase are all interrelated and the whole picture requires further
investigation. After applying Mellin transforms between the theta
and zeta functions and running the logic in reverse, the plan is to unify
the entire three scenarios.

The main component of the Drinfeld Associator is constructed thru
analysis of the KZ equations which arise in Conformal Field Theory.
Numerous presentations of the KZ equations exist. The various formalisms
for representations of solutions of KZ equations in terms of Hypergeometric
functions are well known.

Kontsevich (QA/9904055 [25]) observed the appearance of Drinfeld integrals
and MZV's in several distinct scenarios: Drinfeld's work on the KZ equation
and the associator; Etingof-Kazhdan's quantization of Poisson-Lie algebras;
Tamarkin's proof of the formality theorem; and Kontsevich's quantization
of Poisson manifolds. Kontsevich claims there on page 25 that the appearance
of these various forms of Drinfeld Integrals and MZV's as coefficients are
closely related and unavoidable. In this work, we explain how the various
components are all related.

Utilizing Terasoma's theorem regarding Selberg integrals and MZV's,
and an additional observation regarding Kontsevich's deformation quantization
a complete answer to Kontsevich's remark about the related integrals is indicated.
The coefficients of the terms in f *g in (Cattaneo-Felder QA/9902090 [8]) can be
expressed in terms of combinatorics of Feynman type diagrams. These
combinatorics can in turn be related to the combinatorics used by Terasoma in
his proof of the Selberg-MZV's theorem. Furthermore, this is just one of many
other routes which connect the MZV's and Selberg integrals with Conformal
Field Theory, Associators, and hence the various forms of Drinfeld integrals.
In particular, using the formalism of Beilinson and Drinfeld [2] in their extensive
paper on Chiral Algebras as a starting point, a unifying logical architecture starts
to fall in place.

This analysis is further strengthened by the work of Casati, Falqui, Magri and Pedroni [7] wherein the KP equations are related to the Adler-Gelfand-Dickey brackets on Poisson Manifolds.

1. Mulase's Theorem on Solutions of the KP Equations

    [as stated in the intro above]

2. Le and Murakami's Theorem on Multiple Zeta Values

    [as stated in the intro above]

3. Terasoma's Theorem on Selberg Integrals

    The proof of the main theorem by Terasoma uses Drinfeld's Associator to construct the relation between Selberg integrals and multiple zeta values. While the theorem does not construct an explicit evaluation, it does prove that a Selberg integral has its image in the homogenous ring of MZV's. His proof uses Gelfand generalized hypergeometric functions as an underlying concept. Nevertheless, there still remains substantial applications of the Gelfand apparatus to a more broad view of multiple zeta values, which we do in this paper. Also, the rooted graphs Terasoma constructs to generate the combinatorial Selberg integral appear to be equivalent to the rooted trees used by Connes and Kreimer in their Hopf algebras used to evaluate Feynman diagrams. See Varilly [39].

4. Drinfeld Integrals and Multiple Zeta Values

    See Drinfeld various papers.

5. Etingof Kazhdan Theory:

    In Etingof-Kazhdan IV, the quantum conformal blocks in the WZW model are constructed, from which the quantum KZ connection is constructed. In Etingof-Kazhdan V [14], quantum vertex operator algebras are defined and analyzed. See also Schlichenmaier and Sheinman [36]. Also there are apparently several other different concepts of quantum vertex algebras – from Borcherds (QA/9903038 [3]) and Beilinson-Drinfeld (Chiral Algebras [3]), which require reconciliation, and multiple zeta values might play a role.

6. Drinfeld Integrals and Kontsevich Deformation Quantization

Broadhurst and Kreimer (Vienna [6]) point out the significance of studying Feynman diagrams with only one-loop primitive propagator subdivergences. Aside from numerous references to MZV's and Hopf algebras, they remind us that these Feynman diagrams provide a solution to a universal problem in Hochschild cohomolgy. See Connes-Kreimer (9808042 and 9904044 [10]) for more details. This same Hochschild cohomology is used in the construction of Poisson structure, which is subsequently used in deformation quantization. See Flato-Sternheimer [15].

Furthermore, Cattaneo, Felder and Tomassini (0012228 [9]) give an explicit construction of a star product from the deformation algebra on any Poisson manifold. The deformed algebra of functions is realized as the algebra of horizontal sections of a vector bundle with flat connection. The deformation of transition functions can be expressed in terms of the Ward identities, which correspond to the formality theorem. See also Cattaneo-Felder (9902090 [8]) and Takhtajan [37] on the Ward identities.

Also, the Zero-Curvature-Condition has an equivalent formulation as a Lax Equation. This leads one to suggest there is a greater role in the MZV picture for Lax Equations and yet other techniques of integrable systems. See Palais [33].

7. Beilinson and Drinfeld: Chiral Algebra Theory

A reconciliation between the various versions Etingof Kazhdan [14], Borcherds [3] etc. of vertex algebras would be an appropriate thing to do. In the context a larger view of things where the multiple zeta values play a central role, this is now a more feasible task. Also, recent work of Racinet [34] would come into play here.

8. Calabi-Yau Manifolds

In a theorem of Libgober [28] "the polynomials Q(c,…c) in the Chern classes of certain Calabi-Yau manifolds X are related to the coefficients of the generalized hypergeometric series expansion of the period (holomorphic at maximum degeneracy point) of a mirror of X" {quote from Hoffman}. In (9908065 [20]) Hoffman relates these Chern classes to multiple zeta values.

In an early paper by Klemm et al on Calabi-Yau manifolds (CERN-TH.7369/94 [24]), clues to the theory underlying the connection between CY manifolds and MZV's. The problem to count the "number" of rational curves of arbitrary degree d on the quintic threefold X in $P^4$. The large volume expansion of the correlation function (known as the Yukawa coupling) counts these curves. This function is expressable in a closed form in terms of solutions of a generalized hypergeometric system. This correlation function depends upon the complex structure modulus of X according to the Picard-Fuchs equation. After some technicalities, the Picard-Fuchs ODE is derived from the Dwork-Griffiths-Katz reduction Griffiths [19].

Klemm et. Al. then observes that these ODE are equivalent to the generalized hypergeometric system of Gelfand, Kapranov and Zelevinskii (GKZ1 and GKZ2 [16] and [17]). Also see Hosono, et al. (9511001 [21]). The relation between the Chern classes and the multiple zeta numbers indicates that the period integral which come up require further investigation in this context. See also (9707003 [ ]) and Kontsevich-Zagier (Periods [26]). Also see Richards-Zheng (IAS preprint [35]) for a discussion of Varchenko's closed form expression of the hypergeometric period matrices.

Summary.

This paper presents applications of multiple zeta values to prove a theorem simplifying the evaluation of Hermitian matrix integrals. We extend a result of Mulase wherein he demonstrates that the integrals fall into three categories. Our new theorem proves that all three of Mulase's options are essentially equivalent, and that in all cases, the integral evaluates to a multiple zeta function. We them show how additional considerations of multiple zeta values organize Drinfeld integrals and associators into a well understood theory. We finally show how the Borcherds, Beilinson-Drinfeld, and Etingof-Kazhdan VOA's can be unified and extended to further explain the concepts of multiple zeta values and field theory.

Acknowledgements.

We extend our appreciation and gratitude to Peter Bouwnegt (Univ of Adelaide), Pierre Deligne (IAS), Vladimir Drinfeld (Univ of Chicago), Motohico Mulase (Univ of Calif./Davis) and Donald Richards (IAS/Univ of Va.) for discussions and lectures relevant to this paper.

In particular, Donald conjectured to me in the Spring of 2001 about stronger version of Terasoma's theorem in which Selberg integrals are evaluated explicitly in terms of MZV's. It was this conjecture which led to my deeper investigations into Terasoma's paper. The conjecture is still part of my work in progress. Deligne gave several lectures on Racinet's paper and we had several discussions on MZV's. Mulase explained to me the unpublished details on the proof of his theorem based upon work of Kakei and Goulden-Harer-Jackson, Bouwnegt gave me some advice on the KZ equations and field theory, and Drinfeld gave me a very helpful lecture on associators.

And finally tremendous thanks to Michael Anshel for everything, including his advice and general direction to look towards zeta functions.

David H. Wohl

Hadar HaTorah Yeshiva
824 Eastern Parkway
Brooklyn, NY
11213

Department of Mathematics
Touro College
1602 Avenue J
Brooklyn, NY 11230

DAVIDHW@TOURO.EDU

DHWSYSTEMS@ATT.NET